\documentstyle[amssymb,amsfonts]{amsart}

\def\R{{\hbox{\bf R}}}

\def\P{{\hbox{\bf P}}}
\def\E{{\hbox{\bf E}}}

\font \roman = cmr10 at 10 true pt

\def\tr{{\hbox{\roman tr}}}
\def\diag{{\hbox{\roman diag}}}
\def\dim{{\hbox{\roman dim}}}

\def\dist{{\hbox{\roman dist}}}

\def\Z{{\hbox{\bf Z}}}

\def\eps{\varepsilon}
\newenvironment{proof}{\noindent {\bf Proof} }{\endprf\par}
\def \endprf{\hfill  {\vrule height6pt width6pt depth0pt}\medskip}
\def\emph#1{{\it #1}}
\def\textbf#1{{\bf #1}}

\def\BP{{\mathbf P}}

\def\BBR {{\mathbb R}}

\def\hs{\hfill $\square$}

\def\ep{\epsilon}

\theoremstyle{plain}
  \newtheorem{theorem}[subsection]{Theorem}
  \newtheorem{conjecture}[subsection]{Conjecture}

  \newtheorem{lemma}[subsection]{Lemma}
  \newtheorem{corollary}[subsection]{Corollary}

\theoremstyle{remark}
  \newtheorem{remark}[subsection]{Remark}
  \newtheorem{remarks}[subsection]{Remarks}

\theoremstyle{definition}
  \newtheorem{definition}[subsection]{Definition}

\include{psfig}

\begin{document}

\title{On random $\pm 1$ matrices: Singularity and Determinant }

\author{Terence Tao}
\address{Department of Mathematics, UCLA, Los Angeles CA 90095-1555}
\email{tao@@math.ucla.edu}
\thanks{T. Tao is supported by a grant from the
Packard Foundation.}

\author{Van Vu}
\address{Department of Mathematics, UCSD, La Jolla, CA 92093-0112}
\email{vanvu@@ucsd.edu}

\thanks{V. Vu is an A. Sloan  Fellow and is supported by an NSF Career Grant.}

\begin{abstract} This papers contains two results concerning random  $n \times n$ Bernoulli matrices.
First, we show that with probability tending to one the
determinant has absolute value $\sqrt {n!} \exp(O(\sqrt {n \ln n}
))$. Next, we  prove a new upper bound $.939^n$ on the probability
that the matrix is singular.
\end{abstract}

\maketitle

\section{Introduction}

Let $n$ be a large integer parameter, and let $M_n$ denote a
random $n \times n$ $\pm 1$  matrix (``random'' meaning with
respect to the uniform distribution, i.e., the entries of $M_n$
are i.i.d. Bernoulli random variables).  Throughout the paper, we
assume that $n$ is sufficiently large, whenever needed. We use
$o(1)$ to denote any quantity which goes to zero as $n \to
\infty$, keeping other parameters (such as $\ep$) fixed.

This model of random matrices is of considerable interest in many
areas, including combinatorics, theoretical computer science and
mathematical physics. On the other hand, many basic questions
concerning this model have been open for a long time. In this
paper, we focus on the following two questions:

{\bf \noindent Question 1.} What is the typical value of the
determinant of $ M_n$ ?

\vskip2mm

{\bf \noindent Question 2.} What is the probability that $M_n$ is
singular ?

\vskip2mm

Let us first discuss Question 1.   From Hadamard's inequality, we
have the bound $|\det(M_n)| \leq n^{n/2}$, with equality if and
only if $M_n$ is an Hadamard matrix. However, in general we expect
$|\det(M_n)|$ to be somewhat smaller than $n^{n/2}$.  Indeed, from
the simple estimate\footnote{Indeed, one can prove \eqref{mom} by expanding $\det M_n$ as the sum of $n!$ signs, and observing
that all the covariances vanish.}
\begin{equation}\label{mom}
\E( (\det M_n)^2 ) = n!
\end{equation}
(first observed  by Tur\'an \cite{Tur}), one is led to conjecture
that $|\det(M_n)|$ should be of the order of $\sqrt{n!} =
e^{-n/2+o(n)} n^{n/2}$ with high probability. On the other hand,
even proving that $|\det M_n|$ is typically positive (or
equivalently, that $M_n$ is typically non-singular) is already a
non-trivial task. This task was first done by Koml\'os \cite{Kom1}
(see Theorem \ref{komlos} below).

\vskip2mm

The first main result of this paper shows that with probability
tending to one (as $n$ tends to infinity), the absolute value of
the determinant is very close to $\sqrt {n!}$.

\begin{theorem} \label{determinant}
$$\BP (|\det M_n| \geq  \sqrt {n!}
\exp( -29 n^{1/2} \ln^{1/2} n ) )  = 1-o(1). $$
\end{theorem}

The constant $29$ is generous  but we do not try to optimize it.

\vskip2mm

Note that  from \eqref{mom} and Chebyshev's inequality that
$$\BP (|\det M_n| \leq \omega(n) \sqrt {n!}  ) )   = 1-o(1) $$
for any function $\omega(n)$ which goes to infinity as $n \to
\infty$.  Combining this with the preceding theorem and the
observation that $\det M_n$ is symmetric around the origin, it
follows that for each sign $\pm$, we have the concentration inequality.
$$\det(M_n) = \pm \sqrt {n!} \exp( O(n^{1/2} \ln^{1/2} n ))$$
with probability $1/2 - o(1)$.

Let us now turn to the problem of determining the probability that
$M_n$ is singular. As mentioned above, Komlos showed, in 1967,
that

\begin{theorem}\label{komlos}\cite{Kom1}  $\BP (\det M_n =0) = o(1)$.
\end{theorem}

The task here is to give a precise formula for $o(1)$ in the right
hand side. Since a matrix $M_n$ with two identical (or opposite)
rows or two identical (or opposite) columns is necessarily
singular, it is easy to see that
$$\BP (\det M_n =0) \ge (1 -o(1)) n^2 2^{1-n}.$$
It has often been conjectured (see e.g. \cite{Ol}, \cite{KKS})
that this is the dominant source of singularity.  More precisely,

\begin{conjecture}\label{pconj}
$$ \BP (\det M_n =0) = (1 -o(1)) n^2 2^{1-n}. $$
\end{conjecture}

Prior to this paper, the best partial result concerning this
conjecture is the following, due to Kahn, Koml\'os and Szemer\'edi
\cite{KKS}:

\begin{theorem} \label{kks}\cite{KKS} We have
$\BP (\det M_n =0) \leq (1-\eps+o(1))^n$, where $\eps := .001$.
\end{theorem}

Our second main result is the following improvement of this
theorem:

\begin{theorem} \label{newkks} We have
$\BP (\det M_n =0) \leq (1-\eps+o(1))^n$, where $\eps := .06191
\dots$.
\end{theorem}

This value of $\eps$ is the unique solution in the interval
$(0,1/2)$ to the equation
\begin{equation}
h(\eps) + \frac{\eps}{\log_2 16/15} = 1,\label{fancy-eps}
\end{equation}
where $h$ is the entropy function
\begin{equation}\label{entropy-def}
h(\eps) := \eps \log_2 \frac{1}{\eps} + (1-\eps) \log_2
\frac{1}{1-\eps}.
\end{equation}

We prove Theorem \ref{newkks} in Sections
\ref{part-3}-\ref{halasz}.  Our argument uses several key ideas
from the original proof of Theorem \ref{kks} in \cite{KKS}, but
invoked in a simpler and more direct fashion.  In a sequel to this
paper  \cite{taovu} we shall use more complicated arguments to
improve this value of $\eps$ further, to $\eps = \frac{1}{4}$.

\vskip2mm

This paper is organized as follows.  In Section \ref{dist-sec} we
establish some basic estimates for the distance between a randomly
selected point on the unit cube $\{-1,1\}^n$ and a fixed subspace,
and in Section \ref{dist-sec2} we obtain similar types of estimates in the
case when the subspace is also random.  In Section \ref{part-2} we then apply those estimates to prove
Theorem \ref{determinant}. As a by-product, we also obtain a short
proof of Theorem \ref{komlos}.  We then give the proof of Theorem
\ref{newkks} in Sections \ref{part-3}-\ref{halasz}.

In this paper we shall try to emphasize simplicity. Several
results obtained in these parts can be extended or refined
considerably with  more technical arguments.  In last part of the
paper (Section \ref{part-4}), we will consider some of these
extensions/refinements. In particular, we prove an extension of
Theorem \ref{komlos} and Theorem \ref{determinant} for more
general models of random matrices.

\section{The distance between a random vector and a  deterministic
subspace}\label{dist-sec}

Let $X$ be a random vector chosen uniformly at random from
$\{-1,1\}^n$, thus $X = (\ep_1,\ldots,\ep_n)$ where
$\ep_1,\ldots,\ep_n$ are i.i.d. Bernoulli signs.  Let $W$ be a
(deterministic) $d$-dimensional subspace of $\R^n$ for some $0
\leq d < n$.  In this section we collect a number of estimates
concerning the distribution of the distance $\dist(X,W)$ from $X$
to $W$, which we will then combine to prove Theorem \ref{komlos}
and Theorem \ref{determinant}.

We have the crude estimate
$$ 0 \leq \dist(X,W) \leq \dist(X,0) = \sqrt{n};$$
later we shall see that $\dist(X,W)$ is in fact concentrated
around $\sqrt{n-d}$ (see Lemma \ref{Tal}).

We next recall a simple observation of Odlyzko.

\begin{lemma}\label{odd}\cite{Ol}
$\P( \dist(X,W) = 0 ) \leq 2^{d-n}$.
\end{lemma}

\begin{proof}  Since $W$ has dimension $d$ in $\R^n$, there is a set of $d$ coordinates
which determines all other $n-d$ coordinates of an element of $W$.
But the corresponding $n-d$ coordinates of $X$ are distributed
uniformly in $\{-1,1\}^{n-d}$ (thinking of the other $k$
coordinates of $X$ as fixed).  Thus the constraint $\dist(X,W) =
0$ can only be obeyed with probability at most $2^{d-n}$, as
desired.
\end{proof}

For a variant of Lemma \ref{odd} which gives a lower bound on $\dist(X,W)$ with high
probability, see, Lemma \ref{b-lemma}.  Next, we establish that $\dist(X,W)$ concentrates near $\sqrt{n-d}$.

\begin{lemma} \label{Tal}
Let $W$ be a fixed subspace of dimension $1 \le d \le n-4$ and $X$
a random $\pm 1$ vector. Then
\begin{equation}\label{mean}
\E( \dist(X,W)^2 ) = n-d.
\end{equation}
Furthermore, for any $t >0$
\begin{equation}\label{pdw}
\P( |\dist(X, W) - \sqrt{n-d}| \ge t+2 ) \le  4 \exp(  -t^2/16).
\end{equation}
\end{lemma}

\begin{proof} Let $P = (p_{jk})_{1 \leq
j,k \leq n}$ be the $n \times n$ orthogonal projection matrix from
$\R^n$ to $W$. Let $D = \diag(p_{11},\ldots,p_{nn})$ be the
diagonal component of $P$, and let $A := P-D = (a_{jk})_{1 \leq j,
k \leq n}$ be the off-diagonal component of $P$.  Since $P$ is an
orthogonal projection matrix, we see that $A$ is real symmetric
with zero diagonal. If we write $X = (\ep_1,\ldots,\ep_n)$, then
from Pythagoras's theorem we have
\begin{align*}
\dist(X, W)^2 &= |X|^2 - |PX|^2 \\
&= n - \sum_{j=1}^n \sum_{k=1}^n \ep_j \ep_k p_{jk} \\
&= n - \tr(P) - \sum_{j=1}^n \sum_{k=1}^n \ep_j \ep_k a_{jk} \\
&= n - d - \sum_{j=1}^n \sum_{k=1}^n \ep_j \ep_k a_{jk}.
\end{align*}
This already gives \eqref{mean}, since $a_{jk}$ vanishes on the
diagonal. Set $Y=\sum_{j=1}^n \sum_{k=1}^n \ep_j \ep_k a_{jk}$. It
is easy to see that
$$\E (Y^2)= 2\sum_{1 \le j, k \le n} a_{jk}^2 =2\tr(A^2). $$

\noindent Observe that as $P$ is a projection matrix, the
coefficients $p_{jk}$ are bounded in magnitude by 1, and we have
$$ \sum_{j=1}^n \sum_{k=1}^n p_{jk}^2 = \tr(P^2) = \tr(P) = d.$$
On the other hand

$$\sum_{j=1}^n p_{jj} = \tr (P) = d$$

\noindent so by Cauchy-Schwartz

$$\sum_{j=1}^n p_{jj}^2 \ge d^2/n. $$

 \noindent This implies that

$$\tr(A^2)=  \sum_{j=1}^n \sum_{k=1}^n p_{jk}^2 - \sum_{j=1}^n p_{jj}^2 \leq d -d^2/n \le \min \{d, n-d \}.$$

\noindent Consider the event $\dist(X, W) \ge \sqrt {n-d} +2$.
This probability of this event is bounded from above by

$$\P (\dist^2(X, W) \ge (n-d)+4 \sqrt{n-d}) = \P(Y \ge
4\sqrt{n-d}) \le \P(Y^2 \ge 16 (n-d)). $$

\noindent By Markov's inequality

$$\P(Y^2 \ge 16(n-d)) \le \frac{\E(Y^2)}{16(n-d)} \le \frac{2(n-d)}{16(n-d)} = \frac{1}{8},
$$

\noindent which implies that the  median $M$ of $\dist(X, W)$ is
at most  $\sqrt {n-d}+2$. To bound $M$ from below, consider the
event $\dist(X, W) \ge \sqrt {n-d} -2$. By a similar argument, the
probability of this event is at most

$$ \P(Y \le -4\sqrt{n-d} +4) \le \P(Y^2 \ge 16 (n-d) -32 \sqrt{n-d}
+16).
$$

By Markov's inequality, the last probability is at most

$$\frac{2(n-d)}{16(n-d) -32 \sqrt{n-d}
+16} < \frac{1}{2}, $$

\noindent for all $d \le n-4$. Thus, we can conclude that
$|M-\sqrt{n-d}| \le 2$.

 Since $\dist(X, W)$ is a convex function on
$\{-1,1 \}^n$ with Lipschitz coefficient 1, Talagrand's \cite{Tal}
inequality implies that

$$\P(|\dist(X, W)-M| \ge t) \le 4 \exp(-t^2/16), $$

\noindent for any $t >0$. Since $|M-\sqrt{n-d}| \le 2$, Lemma
\ref{Tal} follows. \end{proof}

\begin{remarks}
 One can deduce a concentration result similar to Lemma \ref{Tal}
using the high moment method; there is also a slightly weaker statement that can be
obtained from Bonami's inequality \cite{Bon}.

One can have a similar statement for the case $d=n-3, n-2$ and
$n-1$. In these cases $\sqrt {n-d} <2 $, so the event $\dist (X,W)
\le \sqrt{n-d}-2$ holds with probability zero. So the median $M$
is between 0 and 3. Therefore, in these cases

$$ \P( \dist(X, W) \ge 3+ t ) \le  4 \exp(  -t^2/16). $$
\end{remarks}

\section{The distance between a random vector and a  random
subspace}\label{dist-sec2}

The estimates in the last section are quite accurate when $n-d$ is
sufficiently large, but do not provide much useful information
when $n-d$ is small (e.g., $n-d =2$). For instance, it does not
show that the distance is (with high probability) not zero in this
case. Indeed, there are some exceptional spaces $W$ (e.g. the hyperplane
of points $(x_1,\ldots,x_n)$ with $x_1=x_2$) which capture a very large fraction of
the points in $\{-1,1\}^n$.  However, in our applications $W$ is a subspace spanned by random vectors
and will thus ``typically'' not be of the exceptional form described above, in which the unit normal contains
many zero coordinates.  In such a case we can still recover good lower bounds on $\dist(X,W)$ with high probability.
More precisely, we have

\begin{lemma} \label{distance0} Let $X$ be a random vector in $\{-1,1\}^n$, let $1 \leq d \leq n-1$ and $W$ a space spanned by
$d$ random vectors in $\{-1,1\}^n$, chosen independently of each other and with $X$. Then we have
$$\P \big( \dist(X, W) \leq \frac{1}{4n} \big) = O(1/ \sqrt{\ln n}). $$
\end{lemma}

\begin{remark} \label{remark:distance0} In fact, as $W$ is spanned
by random vectors, we can fix $X$. The above formulation is,
however, more convenient for the proof. \end{remark}

The remainder of this section will be devoted to the proof of Lemma \ref{distance0}.

Let $1 \leq l \leq n$.  We
say that $W$ is \emph{$l$-typical} if any unit vector
$(w_1,\ldots,w_n) \in W^\perp$ has at least $l$ coordinates whose
absolute values are at least $\frac{1}{2n}$.  In order to prove
Lemma \ref{distance0}, we need the following

\begin{lemma}[$\dist(X,W)$ is large for typical $W$]\label{off}  Let $W$ be a (deterministic) subspace which is
$l$-typical for some $1 \leq l \leq n$.  Then
$$ \P( \dist(X,W) \leq \frac{1}{4n} ) \leq O( \frac{1}{\sqrt{l}} ).$$
\end{lemma}

\begin{proof}  By hypothesis and symmetry, we may assume without loss of generality
that there is a unit normal $(w_1,\ldots,w_n) \in W^\perp$ such
that $|w_1|, \ldots, |w_l| \geq \frac{1}{2n}$. We then see that
\begin{align*}
\P( \dist(X,W) \leq \frac{1}{4n} ) &= \P( |\ep_1 w_1 + \ldots + \ep_n w_n| \leq \frac{1}{4n} ) \\
&\leq \sup_{x \in \R} \P( |\ep_1 w_1 + \ldots + \ep_l w_l - x| \leq \frac{1}{4n} ) \\
&= \sup_{y \in \R} \P( \ep_1 2nw_1 + \ldots + \ep_l 2n w_l \in
[y,y+1] )
\end{align*}
where we have made the substitutions $x := \sum_{l < j \leq n}
\ep_j w_j$ and $y := 2n x - \frac{1}{2}$ respectively.  To
conclude the claim, we invoke the following variant of
the Littlewood-Offord lemma, due to Erd\"os \cite{erdos}:

\begin{lemma} \label{l-o}\cite{erdos} Let $a_1, \dots, a_k$ be  real numbers
with absolute values larger than one. Then for any interval $I$ of
length at most one

$$\P (\sum_{i=1}^k a_i\ep_i \in I) = O(1/ \sqrt k). $$

\end{lemma}

This lemma was proved by Erd\"os using Sperner's lemma. The reader
may want to check Remark \ref{lo-remark} for a different argument.  Lemma \ref{off} immediately follows.
\end{proof}

\noindent We are now ready to prove Lemma \ref{distance0}.

\begin{proof} It suffices to prove the extremal case when $W$ is
spanned by $n-1$ random vectors.  Set $l := \lfloor \frac{\ln
n}{10} \rfloor$. In light of Lemma \ref{off}, we see that it
suffices to show that
\begin{equation}\label{normal}
\P( W \hbox{ is not } l\hbox{-typical} ) = O(1/ \sqrt{\ln n}).
\end{equation}
If $W$ is not $l$-typical, then there exists a unit vector $w$
orthogonal to $W$ with at least $n-l$ coordinates which are less
than $\frac{1}{2n}$ in magnitude. There are ${n \choose n-l} = {n
\choose l}$ such possibilities for these coordinates. Thus by
symmetry we have
$$ \P( W\hbox{ is not } l\hbox{-typical} ) \leq {n \choose l} \P( W\perp w \hbox{ for some } w \in \Omega )$$
where $\Omega$ is the space of all unit vectors $w =
(w_1,\ldots,w_n)$ such that $|w_j| < \frac{1}{2n}$ for all $l < j
\leq n$.

Suppose that $w \in \Omega$ was such that $W\perp w$, then $X_i
\perp w$ for all $1 \leq i \leq n-1$.  Write $X_i = (\ep_{i,1},
\ldots, \ep_{i,n})$, then
$$ \sum_{j=1}^n \ep_{i,j} w_j = 0.$$
Since $\ep_{i,j} = \pm 1$, and $|w_j| < 1/2n$ for $j > l$, we thus
conclude from the triangle inequality that
$$ |\sum_{j=l+1}^n \ep_{i,j} w_j| \leq (n-l) \frac{1}{2n} \leq \frac{1}{2}.$$
On the other hand, we have
\begin{align*}
\sum_{j=1}^l |w_j| &\geq \sum_{j=1}^l |w_j|^2 \\
&= 1 - \sum_{j=l+1}^n |w_j|^2 \\
&\geq 1 - (n-l) (\frac{2}{n})^2 \\
&\geq 1 - \frac{4}{n}.
\end{align*}
Comparing these two inequalities, we see that (for $n > 8$; the cases $n \leq 8$ are of course trivial)
that for each $1 \leq i \leq n-1$, at least one of the
$\ep_{i,j} w_j$ has to be negative.  Thus, if we let
$\ep_1,\ldots,\ep_l$ be signs such that $\ep_j w_j$ is positive
for all $1 \leq j \leq l$, we thus have
$$ (\ep_{i,j})_{1 \leq j \leq l} \neq (\ep_j)_{1 \leq j \leq l} \hbox{ for all } 1 \leq i \leq n-1.$$
Thus we have

\begin{align*} &  \P( W\perp w \hbox{ for some } w \in \Omega ) \\ & \leq
\sum_{\ep_1,\ldots,\ep_l \in \{-1,1\}} \P( (\ep_{i,j})_{1 \leq j
\leq l} \neq (\ep_j)_{1 \leq j \leq l} \hbox{ for all } 1 \leq i
\leq n-1 ).\end{align*}

 Since the $\ep_{i,j}$ are i.i.d. Bernoulli
variables, we have
$$ \P( (\ep_{i,j})_{1 \leq j \leq l} \neq (\ep_j)_{1 \leq j \leq l} \hbox{ for all } 1 \leq i \leq n-1 )
= (1 - 2^{-l})^{n-1}.$$ Putting this all together, we obtain
\begin{align*}
\P( W\hbox{ is not } l-\hbox{typical} ) &\leq {n \choose l} 2^l (1 - 2^{-l})^{n-1} \\
&\leq n^{l+1} 2^l e^{-2^l (n-1)},
\end{align*}
and \eqref{normal} follows by choice of $l$.  This proves Lemma
\ref{distance0}.
\end{proof}

As a consequence of this lemma, we  derive a short proof of
Theorem \ref{komlos}.  Let $X_1, \dots, X_n$ be the row vectors of
$M_n$ and $W_j$ be the subspace spanned by $X_1, \dots, X_j$.
Observe that if $M_n$ is singular, then $X_1,\ldots,X_n$ are
linearly dependent, and thus we have $\dist(X_{j+1},W_j) = 0$ for
some $1 \leq j \leq n-1$. Thus we have
\begin{align*}
\P( \det(M_n) = 0 ) &\leq \sum_{j=1}^{n-1} \P (\dist(X_{j+1}, W_j) = 0) \\
&= \sum_{j=1}^{n-1} \P( \dist(X,W_j) = 0 ).
\end{align*}

From Lemma \ref{odd} we have $\P( \dist(X,W_j) = 0 ) \leq
2^{j-n}$.  Since $\P(\dist(X,W_j) = 0)$ is clearly monotone
increasing in $j$, we obtain the inequality
$$
\P( \det(M_n) = 0 ) \leq 2^{-k} + k \P( \dist(X,W_{n-1}) = 0 )$$
for any $1 \leq k < n$. By the lemma just proved, $\P(
\dist(X,W_{n-1}) = 0 )= O(1/\sqrt {\ln n})$. By choosing $k=
\ln^{1/4} n$

$$2^{-k} + O( k/ \sqrt {\ln n} ) = o(1) $$

\noindent completing the proof.

\section{Proof of Theorem \ref{determinant}}\label{part-2}

For an $n \times n$ matrix $A$, $|\det A|$ is the volume of the
parallelepiped spanned by the row vectors of $A$. If one instead
expresses this volume in terms of base times height, we obtain the
factorization
$$ |\det(M_n)| = \prod_{0 \leq j \leq n-1} \dist(X_{j+1}, W_j).$$
To estimate this quantity, we shall simply control each of the factors
$\dist(X_{j+1}, W_j)$ separately, using the estimates obtained in the
previous two sections.

\noindent We may assume $n$ is large.  Set $d_0= n- \ln^{1/4} n$.
For $1 \le j \le d_0$

$$\gamma_j := 7 \sqrt {\frac{\ln (n-j)}{n-j} } . $$

It is trivial that all $\gamma _j $ are bounded from above by
$1/2$ if $n$ is sufficiently large. Consider  $1\le j \le d_0$.
Assuming that  $W_j$ has dimension $j$, by Lemma \ref{Tal} we have
that the probability that the distance

$$\dist(X_{j+1}, W_j) \le  (1-\gamma_j) \sqrt {n-j} $$

\noindent is at most

$$4 \exp( -\gamma_j^2 (n-j)/16) = 4 \exp (-\frac{49}{16} \ln (n-j)) \le
(n-j)^{-2} $$

\noindent provided that $n-j$ is sufficiently large. This implies
that with probability at least

$$1- \sum_{j=1}^{d_0}  (n-j)^{-2}  =1-o(1)$$   the distance
$\dist(X_{j+1}, W_j)$ is at least $(1-\gamma_j)\sqrt {n-j}$, for
every $1 \le j \le d_0$. (Notice that if $\dist(X_{j+1}, W_j ) >0$
then $W_{j+1}$ has full dimension $j+1$.)

For $d_0 < j \leq n-1$, we are going to use  Lemma \ref{distance0}
to estimate the distances. By this lemma, we have that with
probability at least

$$ 1 - \sum_{d_0 < j \leq n}  O( \frac{1}{\sqrt{\ln n}} )  =1-o(1) $$

the distance $ \dist(X_{j+1}, W_j)$ is at least $ \frac{1}{4n}$
for every $ d_0< j \leq n-1$. (In fact, the bound holds for all $1
\le j \le n-1$.)

 Combining the two  estimates on distances, we see
that with probability $1 - o(1)$,

$$ \prod_{0 \leq j \leq n-1} \dist(X_{j+1}, W_j)  \geq
  \frac{ \sqrt {n!}} {\sqrt
{(n-d_0)!}}
 (\frac{1}{4n})^{n-d_0} \prod_{j=0}^{d_0} (1-\gamma_j).  $$

Since $n-d_0 =o(\ln n)$, the error term $\frac{1} {\sqrt
{(n-d_0)!}} (\frac{1}{4n})^{n-d_0}$ is only $\exp(- o(\ln ^2n ))$.
The main error term comes from the product $\prod_{j=0}^{d_0}
(1-\gamma_j)$. By,  the definition of $\gamma_j$ and the fact that
all $\gamma_j$  are less than $1/2$, we have

$$\prod_{j=1}^{d_0} (1-\gamma_j) \ge \exp(-2 \sum_{j=1}^{d_0}
\gamma_j ) \ge \exp( -14 \sum_{j=1}^{d_0} \sqrt{ \frac{\ln
(n-j)}{n-j}} ) .
$$

\noindent We use a  rough estimate that

$$ \sum_{j=1}^{d_0} \sqrt{ \frac{\ln (n-j)}{n-j} } \le \sqrt {\ln n}
\int_{0}^{n} x^{-1/2} dx = 2 \sqrt{ n \ln n}. $$

\noindent Putting these together, we obtain, with probability
$1-o(1)$, that

\begin{align*}  \prod_{0 \leq j \leq n-1} \dist(X_{j+1}, W_j) &\ge \sqrt {n!}
\exp( -28 n^{1/2} \ln^{1/2} n  + o(\ln^2 n) ) \\&\ge \sqrt {n!}
\exp( -29 n^{1/2} \ln^{1/2} n  )  \end{align*}  proving the
theorem.\hs

\section{Proof of Theorem \ref{newkks}}\label{part-3}

In this section, we denote $N := 2^n$. Our
goal is to prove that $\BP (\det M_n =0) \le N^{-(1+o(1))\eps}$,
where $\eps$ is as in Theorem \ref{newkks}.

 Notice that if $M_n$ is singular, then $X_1,
\ldots, X_n $ span a proper subspace $V$ of $\BBR^n$. The first
(fairly simple) observation is that we can restrict to the case
$V$ is a hyperplane, thanks to the following lemma:

\begin{lemma}\label{hyper}\cite{KKS}  We have
$$ \P( X_1, \ldots, X_n \hbox{ linearly dependent} ) \leq N^{o(1)} \P(
X_1,\ldots,X_n \hbox{ span a hyperplane}).$$
\end{lemma}

\begin{remark} One can replace $N^{o(1)}$ by $1+o(1)$, but this
refinement has no significance in the current
situation.\end{remark}

\begin{proof} If $X_1, \ldots, X_n$
are linearly dependent, then there must exist $0 \leq d \leq
{n-1}$ such that $X_1, \ldots, X_{d+1}$ span a space of dimension
exactly $d$. Since the number of possible $d$ is at most $n =
N^{o(1)}$, it thus suffices to show that
\begin{align*} & \P( X_1, \ldots, X_{d+1} \hbox{ span a space of dimension
exactly } d ) \\ &\leq \hbox{const} \times \P( X_1, \ldots, X_n
\hbox{ span a hyperplane} )
\end{align*}
for each fixed $d$.  However, from Lemma \ref{odd} we see that
\begin{align*} & \P( X_1, \ldots, X_{d+2} \hbox{ span a space of
dimension exactly } d+1 \\
&\quad\quad| X_1, \ldots, X_{d+1} \hbox{ span a space of dimension
exactly } d ) \geq 1 - 2^{d-n},
\end{align*}
and so the claim follows from $n-d-1$ applications of Bayes'
identity.
\end{proof}

In view of this lemma, it suffices to show

$$ \sum_{V, V \hbox{ hyperplane}}   \P(
X_1,\ldots,X_n \,\, \hbox{span} \,\, V) \le N^{-\eps  + o(1)}. $$

Clearly, we may restrict our attention to those hyperplanes $V$
which are spanned by their intersection with   $\{-1,1\}^n$. Let
us call such hyperplanes \emph{non-trivial}.  Furthermore, we call
a hyperplane $H$ \emph{degenerate} if there is a vector $v$
orthogonal to $H$ and at most $\log\log n$ coordinates of $v$ are non-zero.

Fix a hyperplane $V$.  Clearly we have

\begin{equation} \label{estimate1} \P(X_1,\ldots,X_n \,\,\hbox{span} \,\,V) \le \P (X_1,\ldots,X_n \in
V) = \P( X \in V )^n. \end{equation}

The contribution of the degenerate hyperplanes is negligible,
thanks to the following easy lemma (cf. the proof of
\eqref{normal}):

\begin{lemma}\label{deg-crude}  The number of degenerate non-trivial
hyperplanes is at most $N^{o(1)}$.
\end{lemma}

\begin{proof}
If $V$ is degenerate, then there is an integer normal vector $v =
(v_1,\ldots,v_n)$ with at most $\log\log n$ non-zero entries.
There are $\sum_{k \leq \log \log n} {n \choose k} \leq \log \log
n n^{\log\log n} \leq N^{o(1)}$ possible places for the non-zero
entries. By relabeling if necessary we may assume that it is
$v_{1}, \ldots, v_{k}$ which are non-zero for some $1 \leq k \leq
\log \log n$.  Let $\pi: \{-1,1\}^n \rightarrow \{-1,1\}^k$ be the
obvious projection map.  Then $V$ is then determined by the
projections $\{ \pi(X_1),\ldots,\pi(X_n) \}$, which are a subset
of $\{-1,1\}^k$.  The number of such subsets is at most $2^{2^k}
\leq 2^{2^{\log \log n}} = N^{o(1)}$, and the claim
follows\footnote{The above estimates were extremely crude.  In
fact, as shown in \cite{KKS},  one can replace $\log \log n$ with
a quantity as high as $n - 3 \log_2 n$ and still achieve the same
result. }.
\end{proof}

By Lemma \ref{odd}, $\P( X \in V )$ is at most $1/2$ for any
hyperplane $V$, so the contribution of the degenerate non-trivial
hyperplanes to $\BP (\det M_n =0)$ is only $N^{-1+o(1)}$.

\vskip2mm

Following \cite{KKS}, it will be useful to specify the magnitude
of $\P( X \in V )$. For each non-trivial hyperplane $V$, define
the \emph{discrete codimension} $d(V)$ of $V$ to be the unique
integer multiple of $1/n$ such that
\begin{equation}\label{codimension-def}
 N^{-\frac{d(V)}{n}-\frac{1}{n^2}} < \P( X \in V) \leq
N^{-\frac{d(V)}{n}}.
\end{equation}
Thus $d(V)$ is large when $V$ contains few elements from $\{-1,1\}^n$, and conversely.

We define by $\Omega_d$ the set of all non-degenerate, non-trivial
hyperplanes with discrete codimension $d$. It is simple to see
that $ 1 \leq d(V) \leq n$ for all non-trivial $V$. In particular,
there are at most $O(n^2) = N^{o(1)}$ possible values of $d$, so
to prove our theorem it suffices to prove that

\begin{equation}\label{d-dim}
 \sum_{V \in \Omega_d} \P( X_1,\ldots,X_n \hbox{ span } V)  \leq N^{-\eps
+ o(1)}
\end{equation}
for all $1 \leq d \leq n$.  (Our errors $o(1)$ shall decay to zero as $n \to \infty$
uniformly in the choice of $d$.)

We first handle the (simpler) case when $d$ is large. Note that if
$X_1,\ldots,X_n$ span $V$, then some subset of $n-1$ vectors
already spans $V$. By symmetry, we have

\begin{align*}
\sum_ {V \in \Omega_d} \P( X_1,\ldots,X_n \hbox{ span } V) &\leq n
\sum_{V \in \Omega_d}  \P( X_1,\ldots,X_{n-1}  \hbox{ span } V)
\P(X_n \in V)
\\ &\le nN^{-\frac{d}{n}}\sum_{V \in \Omega_d}  \P(
X_1,\ldots,X_{n-1}  \hbox{ span } V) \\ &\le  nN^{-\frac{d}{n}} =
N^{-\frac{d}{n} + o(1)}
\end{align*}

This disposes of the case when $d \geq (\eps - o(1)) n$.  Thus to
prove Theorem \ref{newkks} it will now suffice to prove

\begin{lemma} \label{mainbound}  If $d$ is any integer multiple of $1/n$ such that
\begin{equation}\label{d-bound}
1 \leq d \leq (\eps - o(1))n
\end{equation}
then we have
$$\sum_{V \in \Omega_d}
\P( X_1,\ldots,X_n \hbox{ span } V) \le  N^{-\eps+o(1)}.$$
\end{lemma}

This is the objective of the next section.

\section{Proof of Lemma \ref{mainbound}}\label{part-3a}

The key idea in \cite{KKS} is to find a new kind of random vectors
which are more concentrated on hyperplanes in $\Omega_d$ (with
small $d$) than $(\pm 1)$ vectors. Roughly speaking, if we can
find a random vector $Y$ such that for any $V \in \Omega_d$

$$\P(X \in V) \le c\P(Y \in V) $$
for some $0 < c < 1$, then, intuitively,  one may expect that

\begin{equation} \label{crude} \P(X_1,\ldots,X_n \hbox{ span } V) \le c^n
\P( Y_1,\ldots,Y_n \hbox{ span } V) \end{equation}

\noindent where $X_i$ and $Y_i$ are independent samples of $X$ and
$Y$, respectively.   Since $Y_1,\ldots,Y_n$ can only span at most one hyperplane $V$,
one can then hope to conclude a bound of $O(c^n)$ for the probability that $X_1,\ldots,X_n$ span a hyperplane.

While (\ref{crude}) may be too optimistic (because the samples of
$Y$ on $V$ may be too linearly dependent), it has turned out that
something little bit weaker  can be obtained, with a proper
definition of $Y$.  We next present this important definition.

\begin{definition}
For any $0 \leq \mu \leq 1$, let $\eta^{(\mu)} \in \{-1,0,1\}$ be
a random variable which takes $+1$ or $-1$ with probabilities
$\frac{\mu}{2}$, and $0$ with probability $1-\mu$.  Let $X^{(\mu)}
\in \{-1,0,1\}^n$ be a random variable of the form $X^{(\mu)} =
(\eta^{(\mu)}_1, \ldots, \eta^{(\mu)}_n)$, where the
$\eta^{(\mu)}_j$ are iid random variables with the same
distribution as $\eta^{(\mu)}$.
\end{definition}

Thus $X^{(1)}$ has the same distribution as $X$, while $X^{(0)}$
is concentrated purely at the origin.  The other random variables $X^{(\mu)}$ have an intermediate
behavior.  We shall work with
$X^{(\mu)}$ for $\mu := 1/16$; this is not the optimal value of
$\mu$ but is the cleanest to work with.  For this value of $\mu$
we have the crucial inequality, following a Fourier-analytic argument of Hal\'asz
\cite{Hal} (see also \cite{KKS}).

\begin{lemma}\label{hl}  Let $V$ be a non-degenerate
non-trivial hyperplane.  Then we have
$$ \P( X \in V ) \leq (\frac{1}{2} + o(1)) \P( X^{(1/16)} \in V ).$$
\end{lemma}

This lemma can be viewed as an assertion that any subspace which contains many points from $\{-1,1\}^n$, must
necessarily contain several further points from $\{-1,0,1\}^n$ (weighted appropriately).
We will prove this lemma in the next section.

\begin{remark} One can obtain similar results for smaller values of
$\mu$ than $1/16$; for instance this was achieved in \cite{KKS}
for the value $\mu := \frac{1}{108} e^{-1/108}$, eventually
resulting in their final gain $\eps :=.001$ in Theorem \ref{kks}.
However the smaller one makes $\mu$, the smaller the final bound
on $\eps$; indeed, most of the improvement in our bounds over
those in \cite{KKS} comes from increasing the value of $\mu$.  One
can increase the $1/16$ parameter somewhat at the expense of
worsening the $\frac{1}{2}$ factor; in fact one can increase
$1/16$ all the way to $1/4$ but at the cost of replacing $1/2$
with $1$.  This shows that $(3/4 + o(1))^n$ is the limit of our
method.  We have actually been able to attain this limit; see
\cite{taovu}.
\end{remark}

Let $V$ be a  hyperplane in $\Omega_d$  for some $d$ obeying the
bound in Lemma \ref{mainbound}.
 Let $\gamma$ denote the
quantity
\begin{equation}\label{gamma-def}
\gamma := \frac{d}{n \log_2 16/15};
\end{equation}
note from \eqref{fancy-eps} and \eqref{d-bound} that $0 < \gamma <
1$. Let $\eps' := \min(\eps,\gamma)$.

Consider the event that the i.i.d random vectors $X_1,\ldots,X_{(1-\gamma) n}, X'_1, \ldots, X'_{(\gamma-\eps') n}$
are linearly independent in $V$ (we omit the rounding which plays
no significant role). One can lower bound the probability of this
event by the probability that all $X_i$ and all  $X'_j$ belong to
$V$, which is

$$\P( X \in V)^{(1-\eps')n} = N^{-(1-\eps') d - o(1)}. $$

Let us replace $X_j$ by $X^{(1/16)}$ for $1 \le j \le (1-\gamma)n$
and consider the event $A_V$ that $
X^{(1/16)}_1,\ldots,X^{(1/16)}_{(1-\gamma) n}, X'_1, \ldots,
X'_{(\gamma-\eps') n}$ are linearly independent in $V$. Using
Lemma \ref{hl}, we are able to give a much better lower bound for
this event:

\begin{equation} \label{boundA} \P(A_V)  \geq N^{(1-\gamma) - (1-\eps') d - o(1)}.
\end{equation}

The critical gain is the term $N^{(1-\gamma)} $. In a sense, this
gain is expected since  $X^{(1/16)}$ is much more concentrated on
$V$ then $X$. We will prove (\ref{boundA}) at the end of the
section. Let us now use it to conclude the proof\footnote{The argument below
is a simplified version of a hypergraph covering argument used in \cite{KKS}.} of Lemma
\ref{mainbound}.

Fix $V \in \Omega_d$. Let us denote by $B_V$ the event that
$X_1,\ldots,X_n \hbox{ span } V$. To prove Lemma \ref{mainbound}, we thus need to show
$$ \sum_{V \in \Omega_d} \P( B_V) \leq N^{-\eps+o(1)}.$$

The idea is to use \eqref{boundA} to ``replace'' some of the $X_1,\ldots,X_n$ with
the random variables $X^{(1/16)}_1,\ldots,X^{(1/16)}_{(1-\gamma)n}$, which are more concentrated on $V$, in order
to obtain an exponential type gain.  Since $A_V$ and $B_V$ are
independent, we have, by (\ref{boundA}) that

$$\P(B_V) =\frac{\P(A_V \wedge B_V)}{\P(A_V)} \le N^{- (1-\gamma) +
(1-\eps')d + o(1)} \P(A_V \wedge B_V). $$

Consider a set
$$X^{(1/16)}_1,\ldots,X^{(1/16)}_{(1-\gamma) n}, X'_1, \ldots,
X'_{(\gamma-\eps') n}, X_1,\ldots,X_n$$ of vectors satisfying $A_V
\wedge B_V$. Then there exists $\eps' n - 1$ vectors $X_{j_1},
\ldots, X_{j_{\eps' n-1}}$ inside $X_1,\ldots,X_n$ which, together
with $X^{(1/16)}_1,\ldots,X^{(1/16)}_{(1-\gamma) n}, X'_1, \ldots,
X'_{(\gamma-\eps') n}$, span $V$. Since the  number of possible
indices $(j_1,\ldots, j_{\eps' n-1})$ is ${n \choose \eps' n - 1}
= N^{h(\eps') + o(1)}$,   by conceding a factor of $N^{h(\eps') +
o(1)}$, we can assume that $j_i =i$ for all relevant $i$. Let
$C_V$ be the event that
$X^{(1/16)}_1,\ldots,X^{(1/16)}_{(1-\gamma) n}, X'_1, \ldots,
X'_{(\gamma-\eps') n}, X_1, \ldots, X_{\eps' n-1} \hbox{ span }
V$. Then we have

$$\P(B_V) \le  N^{-(1-\gamma) +
(1-\eps')d + h(\eps') + o(1)} \P\Big(C_V \wedge (X_{\eps'
n},\ldots,X_n \hbox{ in } V ) \Big). $$

\noindent On the other hand,  $C_V$ and the event $(X_{\eps'
n},\ldots,X_n \hbox{ in } V )$ are independent, so

$$\P\Big(C_V \wedge (X_{\eps' n},\ldots,X_n \hbox{ in } V ) \Big)
=\P(C_V) \P( X \in V )^{(1-\ep')n+1}. $$

\noindent Putting the last two estimates together we obtain

\begin{align*}  \P(B_V) &\leq N^{-(1-\gamma) + (1-\eps')d + h(\eps') + o(1)}
N^{-((1-\eps') n + 1)d/n} \P(C_V) \\ & =N^{-(1-\gamma) +h(\ep') -
\ep+o(1)} \P(C_V). \end{align*}

\noindent Since any set of vectors can only span a single space
$V$, we have $\sum_{V \in \Omega_d} \P(C_V) \le 1$. Thus, by
summing over $\Omega_d$, we have

$$\sum_{V \in \Omega_d} \P(B_V) \le N^{-(1-\gamma) +h(\ep') -
\ep+o(1)}. $$

\noindent  We can rewrite the right hand side  using
\eqref{gamma-def} as $ N^{h(\eps') + \frac{d}{n} (\frac{1}{\log_2
16/15} - 1) - 1 + o(1)}. $ Since $\frac{1}{\log_2 16/15} - 1 > 0$,
$d/n \le \ep$, and $h$ is monotone in the interval $0 < \eps' \leq
\eps < 1/2$ we obtain

$$ \sum_{V \in \Omega_d} \P( B_V)  \leq N^{h(\eps) + \eps (\frac{1}{\log_2 16/15} -
1) - 1 + o(1)}.
$$
and the claim follows from the definition of $\ep$ in
\eqref{fancy-eps}. \hs

\noindent In the rest of this section, we prove (\ref{boundA}).
The proof of Lemma \ref{hl}, which uses entirely different
arguments (based on Fourier analysis), will be presented in the next section.

To prove (\ref{boundA}), first notice that the right hand side is
the probability of the event $A'_V$ that $ X^{(1/16)}_1, \ldots,
X^{(1/16)}_{(1-\gamma)n}, X'_1, \ldots, X'_{(\gamma - \eps)n}$
belong to $ V$. Thus, by Bayes' identity it is sufficient to show
that
$$\P(A_V|A'_V) = N^{o(1)}. $$

This is an estimate similar to that in Lemma \ref{hyper}, and we prove it by similar
arguments.
From \eqref{codimension-def} we have
\begin{equation}\label{V-1}
 \P( X \in V ) = (1 + O(1/n)) 2^{-d}
 \end{equation}
and hence by Lemma \ref{hl}
\begin{equation}\label{V-16}
 \P( X^{(1/16)} \in V ) \geq (2 + O(1/n)) 2^{-d}.
 \end{equation}
On the other hand, by a trivial modification of the proof of Lemma
\ref{odd} we have
$$ \P( X^{(1/16)} \in W ) \leq (15/16)^{n - \dim(W)}$$
for any subspace $W$.  By Bayes' identity we thus have the
conditional probability bound
$$ \P( X^{(1/16)} \in W | X^{(1/16)} \in V) \leq (2 + O(1/n)) 2^d
(15/16)^{n - \dim(W)}.$$ This is non-trivial when $\dim(W) \leq
(1-\gamma)n$ thanks to \eqref{gamma-def}.

Let $E_k$ be the event that $X^{(1/16)}_1, \ldots, X^{(1/16)}_k$
are independent. The above estimates imply that

$$\P(E_{k+1} | E_k \wedge A'_V)
 \geq
1 - (2 + O(1/n)) 2^d (15/16)^{n - k}. $$

\noindent for all $0 \leq k \leq (1-\gamma)n$. Applying Bayes'
identity repeatedly (and \eqref{gamma-def}) we thus obtain
$$  \P( E_{(1-\gamma)n} | A'_V ) \geq
N^{-o(1)}.$$

If $\gamma \leq \eps$ then we are now done, so suppose $\gamma >
\eps$ (so that $\eps' = \eps$). From Lemma \ref{odd} we have
$$ \P( X \in W ) \leq (1/2)^{n - \dim(W)}$$
for any subspace $W$, and hence by \eqref{V-1}
$$ \P( X \in W | X \in V) \leq (1 + O(1/n)) 2^d (1/2)^{n-\dim(W)}.$$

Let us assume $E_{(1-\gamma n)}$ and denote by $W$ the  $(1-\gamma
n)$-dimensional subspace spanned by $X^{(1/16)}_1, \ldots,
X^{(1/16)}_{(1-\gamma)n}$. Let $U_k$ denote the event that $ X'_1,
\ldots, X'_{k}, W \hbox{ are independent}$. We have

$$p_k= \P(U_{k+1}|U_k \wedge A'_V) \ge  1 - (1 + O(1/n)) 2^d (1/2)^{n-k-(1-\gamma)n} \ge
 1 - \frac{1}{100} 2^{(k+\eps-\gamma)n} $$

\noindent  for all $0 \leq k < (\gamma - \eps)n$, thanks to
\eqref{d-bound}. Thus by Bayes' identity we obtain

$$\P(A_V|A'_V) \ge N^{o(1)}  \prod_{0 \leq k < (\gamma - \eps)n} p_k =
N^{o(1)} $$  as desired. \hs

\section{Hal\'asz-type arguments}\label{halasz}

We now prove Lemma \ref{hl}.  The first step is to use Fourier
analysis as in \cite{Hal} to obtain usable formulae for $\P(X \in V)$ and $\P(
X^{(\mu)} \in V )$.  Let $v \in \Z^n \backslash \{0\}$ be an
normal vector to $V$ with integer coefficients (such a vector
exists since $V$ is spanned by the integer points $V \cap
\{-1,1\}^n$). By hypothesis, at least $\log \log n$ of the coordinates of
$v$ are non-zero.

We first observe that the probability $\P( X^{(\mu)} \in V )$ can
be computed using the Fourier transform:
\begin{align*}
\P(X^{(\mu)} \in V ) &= \P(X^{(\mu)} \cdot v = 0)\\
 &= \E( \int_{0}^1 e^{2\pi i \xi X^{(\mu)} \cdot v}\ d\xi ) \\
&= \int_0^1 \E( e^{2\pi i \xi \sum_{j=1}^n \epsilon^{(\mu)}_j
v_j})
\ d\xi\\
&= \int_0^1 \prod_{j=1}^n ((1-\mu) + \mu \cos(2\pi \xi v_j))\
d\xi.
\end{align*}
Applying this with $\mu=1/16$ we obtain
$$ \P(X^{(1/16)} \in V ) = \int_0^1 \prod_{j=1}^n (\frac{15}{16} +
\frac{1}{16} \cos(2\pi \xi v_j))\ d\xi.$$ Applying instead with
$\mu=1$, we obtain
\begin{align*}
\P(X \in V) &= \int_0^1 \prod_{j=1}^n \cos(2\pi \xi v_j)\ d\xi \\
&\leq \int_0^1 \prod_{j=1}^n |\cos(2\pi \xi v_j)|\ d\xi \\
&= \int_0^1 \prod_{j=1}^n |\cos(\pi \xi v_j)|\ d\xi,
\end{align*}
where the latter identity follows from the change of variables
$\xi \mapsto \xi/2$ and noting that $|\cos(\pi \xi v_j)|$ is still
well-defined for $\xi \in [0,1]$.  Thus if we set
\begin{equation}\label{fg-def}
F(\xi) := \prod_{j=1}^n |\cos(\pi \xi v_j)|; \quad G(\xi) :=
\prod_{j=1}^n (\frac{15}{16} + \frac{1}{16} \cos(2\pi \xi v_j)),
\end{equation}
it will now suffice to show that
\begin{equation}\label{G-market}
\int_0^1 F(\xi)\ d\xi \leq (\frac{1}{2}+o(1)) \int_0^1 G(\xi)\
d\xi.
\end{equation}
We now observe three estimates on $F$ and $G$.

\begin{lemma} \label{3estimates} For any $\xi, \xi' \in [0,1]$, we have the pointwise
estimates
\begin{equation}\label{fag}
 F(\xi) \leq G(\xi)^4
 \end{equation}
and
\begin{equation}\label{ffg}
F(\xi) F(\xi') \leq G(\xi + \xi')^2
\end{equation}
and the crude integral estimate
\begin{equation}\label{G-int}
 \int_0^1 G(\xi)\ d\xi \leq o(1).
\end{equation}
\end{lemma}

Of course, all operations on $\xi$ and $\xi'$ such as  $(\xi +
\xi')$ in (\ref{ffg}) are considered  modulo $1$.

\vskip2mm

{\bf \noindent Proof of Lemma \ref{3estimates}.}   We first prove
\eqref{fag}.  From \eqref{fg-def} it will suffice to prove the
pointwise inequality
$$ |\cos \theta| \leq [ \frac{15}{16} + \frac{1}{16} \cos 2\theta
]^4$$ for all $\theta \in \R$.  Writing $\cos 2\theta = 1-2x$ for
some $0 \leq x \leq 1$, then $|\cos \theta| = (1-x)^{1/2}$ and the
inequality becomes
$$ (1-x)^{1/2} \leq (1 - x/8)^{4}.$$
Introducing the function $f(x) := \log(\frac{1}{1-x})$, this
inequality is equivalent to
$$ \frac{f(x) - f(0)}{x-0} \geq \frac{f(x/8) - f(0)}{x/8 - 0}$$
but this is immediate from the convexity of $f$.

Now we prove \eqref{ffg}.  It suffices to prove that
$$ |\cos \theta| |\cos \theta'| \leq [ \frac{15}{16} + \frac{1}{16}
\cos(2(\theta+\theta')) ]^2$$ for all $\theta, \theta' \in \R$. As
this inequality is periodic with period $\pi$ in both $\theta$ and
$\theta'$ we may assume that $|\theta|, |\theta'| < \pi/2$ (the
cases when $\theta = \pi/2$ or $\theta' = \pi/2$ being trivial).
Next we observe from the concavity of $\log \cos(\theta)$ in the
interval $(-\pi/2,\pi/2)$ that
$$ \cos \theta \cos \theta'  \leq \cos^2 \frac{\theta + \theta'}{2} =
\frac{1}{2} + \frac{1}{2} \cos(\theta + \theta').$$ Writing
$\cos(\theta + \theta') = 1 - 2x$ for some $0 \leq x \leq 1$, then
$\cos2(\theta+\theta') = 2(1-2x)^2 - 1 = 1 - 8x + 8x^2$, and our
task is now to show that
$$ 1 - x \leq (1 - (x-x^2)/2)^2 = 1 - x + x^2 + (x-x^2)^2/4,$$
but this is clearly true.

Now we prove \eqref{G-int}.  We know that at least $\log \log n$ of the
$v_j$ are non-zero; without loss of generality we may assume that
it is $v_1, \ldots, v_K$ which are non-zero for some $K > \log \log n$.  Then we have by
H\"older's inequality, followed by a rescaling by $v_j$
\begin{align*}
\int_0^1 G(\xi)\ d\xi &\leq \int_0^1 \prod_{j=1}^{K}
(\frac{15}{16} + \frac{1}{16} \cos(2\pi \xi v_j))\ d\xi \\
&\leq \prod_{j=1}^{\log\log n} (\int_0^1 (\frac{15}{16} + \frac{1}{16}
\cos(2\pi \xi v_j))^{\log\log n}\ d\xi)^{1/\log\log n} \\
&= \prod_{j=1}^{K} (\int_0^1 (\frac{15}{16} + \frac{1}{16} \cos(2\pi
\xi))^{K} \ d\xi)^{1/{K}} \\
&= \int_0^1 (\frac{15}{16} + \frac{1}{16} \cos(2\pi \xi))^{K}\ d\xi\\
&= o(1)
\end{align*}
as desired, since $K \geq \log\log n$.
\hs

Now we can quickly conclude the proof of \eqref{G-market}.  From
\eqref{ffg} we have the sumset inclusion
$$
\{ \xi \in [0,1]: F(\xi) > \alpha \} + \{ \xi \in [0,1]: F(\xi) >
\alpha \} \subseteq \{ \xi \in [0,1]: G(\xi) > \alpha \}$$ for any
$\alpha > 0$.  Taking measures of both sides and applying the
Mann-Kneser-Macbeath ``$\alpha+\beta$ inequality'' $|A+B| \geq
\min(|A|+|B|, 1)$ (see \cite{macbeath}), we obtain
$$ \min( 2|\{ \xi \in [0,1]: F(\xi) > \alpha \}|, 1)
\leq  |\{ \xi \in [0,1]: G(\xi) > \alpha \}|.$$ But from
\eqref{G-int} we see that $|\{ \xi \in [0,1]: G(\xi) > \alpha \}|$
is strictly less than 1 if $\alpha > o(1)$.  Thus we conclude that
$$ |\{ \xi \in [0,1]: F(\xi) > \alpha \}| \leq \frac{1}{2} |\{ \xi \in
[0,1]: G(\xi) > \alpha \}|$$ when $\alpha > o(1)$.  Integrating
this in $\alpha$, we obtain
$$
\int_{[0,1]: F(\xi) > o(1)} F(\xi)\ d\xi \leq \frac{1}{2} \int_0^1
G(\xi)\ d\xi.$$ On the other hand, from \eqref{fag} we see that
when $F(\xi) \leq o(1)$, then $F(\xi) = o( F(\xi)^{1/4}) \leq
G(\xi)$, and thus
$$
\int_{[0,1]: F(\xi) \leq o(1)} F(\xi)\ d\xi \leq o(1) \int_0^1 G\
d\xi.$$ Adding these two inequalities we obtain \eqref{G-market}
as desired.  This proves Lemma \ref{hl}. \hs

\begin{remark}\label{lo-remark} A similar Fourier-analytic argument can be used to prove
Lemma \ref{l-o}.  To see this, we first recall Ess\'een's
concentration inequality \cite{esseen}
$$ \P( X \in I ) \leq C \int_{|t| \leq 1} |\E( e^{i t X} )|\ dt$$
for any random variable $X$ and any interval $I$ of length at most
1.  Thus to prove Lemma \ref{l-o} it would suffice to show that
$$ \int_{|t| \leq 1} |\E( \exp( it \sum_{j=1}^k a_j \ep_j )| )\ dt = O( 1 / \sqrt{k} ).$$
But by the independence of the $\ep_j$, we have
$$ |\E( \exp( it \sum_{j=1}^k a_j \ep_j ) )|  =
\prod_{j=1}^k |\E( e^{i t a_j \ep_j} )| = |\prod_{j=1}^k \cos(t
a_j)|$$ and hence by H\"older's inequality
$$ \int_{|t| \leq 1}| \E( \exp( it \sum_{j=1}^k a_j \ep_j ))|\ dt
\leq \prod_{j=1}^k (\int_{|t| \leq 1} |\cos(t a_j)|^k\
dt)^{1/k}.$$ But since each $a_j$ has magnitude at least 1, it is
easy to check that $\int_{|t| \leq 1} |\cos(t a_j)|^k\ dt = O( 1 /
\sqrt{k} )$, and the claim follows.
\end{remark}

\section{Extensions and Refinements}\label{part-4}

\subsection{Singularity of more general random matrices}
In \cite{Kom2}, Koml\'os extended Theorem \ref{komlos} by showing
that the singularity probability is still $o(1)$ for a random
matrix whose entries are i.i.d. random variables with
non-degenerate distribution. By slightly modifying our proof of
Theorem \ref{komlos}, we are able to prove a different extension.

We say that a random variable $\xi$ has $(c, \rho)$-property if

$$\min \{ \P (\xi \ge c), \P (\xi \le -c) \} \ge \rho. $$

Let $\xi_{ij}$, $1 \le i, j \le n$ be independent random
variables. Assume that there are positive constants $c$ and $\rho$
(not depending on $n$) such that for all $1 \le i, j \le n$,
$\xi_{ij}$ has $(c,\rho)$-property. The new feature here is that
we {\it do not} require $\xi_{ij}$ be identical.

\begin{theorem} \label{komlos2} Let $\xi_{ij}$, $1 \le i, j \le n$
be as above. Let $M_n$ be the random matrix with entries
$\xi_{ij}$. Then

$$\P(\det M_n =0) = o(1). $$

\end{theorem}

We only  sketch the proof, which follows the proof of Theorem
\ref{komlos}  very closely and uses the same notation: $X_1,
\dots, X_n$ are the row vectors of $M_n$ and $W_j$ is the subspace
spanned by $X_1, \dots, X_j$. We will show

\begin{equation} \label{sumWjgeneral}
\sum_{j=1}^{n-1} \P (X_{j+1} \in W_j) =o(1). \end{equation}

 \noindent This estimate  is a consequence of the following two
 lemmas, which are generalization of Lemmas \ref{odd} and
 \ref{distance0}.

\begin{lemma} \label{oddgeneral} Let $W$ be a $k$ dimensional subspace of $\BBR^n$. Then
for any $1\le j \le n$

$$\P( X_j \in W) \le (1- \rho)^{n-k}. $$
 \end{lemma}

\begin{lemma} \label{lastgeneral} For any $n/2\le j \le n$

$$\P (X_j \in W_{j-1}) = O(1/\sqrt {\ln n} ). $$

\end{lemma}

The proof of Lemma \ref{oddgeneral} is the same as that of Lemma
\ref{odd}. The only information we need is that for any fixed
number $x$ and any plausible $i,j$, $\P(\xi_{ij} =x) \le 1- \rho$.

To prove Lemma \ref{lastgeneral}, let us consider the case $j=n$
(the proof is the same for other cases). We need to modify the
definition of universality as follows.

 We  call a subset $V $
of $n$-dimensional vectors $k$-{\it universal} if for any set of
$k$ indices $1 \le i_1 < i_2 < \dots < i_k \le n$ and any sign
sequence $\ep_1, \dots, \ep_k$, one can find a vector $v \in V$,
such that  the $i_j$ coordinate of $v$ has sign $\ep_j$ and
absolute value  at least $c$.

\vskip2mm

In what follows, we set $l= \ln n/10$. We first show $X_1, \dots,
X_n$ is very likely to be $l$-universal. (Notice that the $X_j$
have different distribution.)

\vskip2mm

\begin{lemma} \label{universalgeneral}  With probability $1-o(1/n)$, $X_1, \dots, X_n$
 is $l$-universal. \end{lemma}

 {\bf \noindent Proof of Lemma \ref{universalgeneral}.} Fix a set
 of indices and a sequence of signs. For any $1\le j \le n$, the
 probability that $X_j$ fails is at most $1- \rho^{l}$. The rest
 of the proof is the same. \hs

 \noindent It follows that

 \begin{corollary}\label{normalgeneral}  Let $H$ be
a subspace  spanned by $n-1$ random vectors. Then with probability
$1-o(1/n)$, any unit vector perpendicular to  $H$ has at least
$l+1$ coordinates whose absolute values are at least
$\frac{1}{Kn}$, where $K$ is a constant depending on $c$.
\end{corollary}

The last ingredient is the following generalization of Lemma
\ref{l-o}.

\begin{lemma} \label{l-ogeneral} Let $a_1, \dots, a_k$ be  real numbers
with absolute values larger than one and $\ep_1, \dots, \ep_k$ be
independent random variables satisfying the $(c, \rho)$-property.
Then for any interval $I$ of length one

$$\P (\sum_{i=1}^k a_i\ep_i \in I) = O(1/ \sqrt k). $$

\end{lemma}

Theorem \ref{komlos2} follows from Corollary \ref{normalgeneral}
and Lemma \ref{l-ogeneral}. To conclude, let us remark that
statements more accurate than   Lemma \ref{l-ogeneral} are known
(see e.g. \cite{Hal}). However, this lemma can be proved using an
argument similar to the one in Remark \ref{lo-remark}.

\vskip2mm

\subsection{Determinants of more general random matrices}

Let $\xi_{ij}$, $1\le i,j \le n$, be a set of independent (but not
necessarily i.i.d.) r.v's with the following two properties:

\begin{itemize}

\item Each $\xi_{ij}$ has mean zero and variance one.

\item There is a constant $K$ that $|\xi_{ij}| \le K$ with
probability one.

\end{itemize}

These two properties imply the following property

\begin{itemize}
\item There are constants $\delta >0$ and $\delta' >0$ such that
for any interval $I$ of length $2\delta$, $\P(\xi_{ij } \in I) \le
1-\delta'$ for all $1\le i,j \le n$.

\end{itemize}

\begin{theorem} \label{determinantgeneral}  
Consider the random matrix $M_n $ with
entries $\xi_{ij}$ as above.
Let $\ep$ be an arbitrary  positive constant.
With probability $1-o(1)$,

$$ |\det M_n| \geq  \sqrt {n!} \exp( -n^{1/2+\ep}  ). $$

\end{theorem}

Notice that
 Lemma \ref{odd} holds for this
model of random matrices, since the last property of $\xi_{ij}$
implies that $\xi_{ij}$ has $(c, \rho)$ property.

Next, consider  Lemma \ref{Tal}. Consider a row vector, say,
$X=(\xi_{i1}, \dots, \xi_{in})$ and a fixed subspace $W$ of
dimension $d$. Again, we have (with the same notation as in
Section 2)

$$ \dist(X, W)^2 = |X|^2 - |PX|^2 = |X|^2 - \sum_{j=1}^n \sum_{k=1}^n \xi_{ij}
\xi_{ik} p_{jk}. $$

\noindent However, it is  no longer the case that the last formula equals

$$ n - d - \sum_{j=1}^n \sum_{k=1}^n \xi_{ij} \xi_{ik} a_{jk} $$

\noindent since $\xi_{ij}$ are not Bernoulli random variables. On
the other hand, we can have something similar with an extra error
term. It is easy to show, using Chernoff's bound, that

$$|X|^2 =\sum_{j=1}^n \xi_{ij}^2 \ge n - \frac{C}{2} n^{1/2} \ln n $$
holds with probability at least $1- 1/2n^2$, for some sufficiently
large $C$. Similarly,

$$\sum_{j=1}^n \xi_{ij}^2 p_{jj} \le d- \frac{C}{2} n^{1/2} \ln n
$$

\noindent holds with probability at least $1- 1/2n^2$. (The use of
Chernoff's bound requires of random variables be bounded. One can
of course, use some other method to remove this assumption.)

The probability $1/n^2$ is negligible. Moreover, we can apply
Talagrand's inequality the same way as before. However, because of
the new error term $C n^{1/2} \ln n$, we cannot set $d_0= n -
\ln^{1/4} n$, but have to stop at $ n- Cn^{1/2} \ln n$. In order
to handle  the cases when $\ln^{1/4} n \le n-d \le C n^{1/2} \ln
n$, we need the following lemma, due to Bourgain (private
conversation), which can be seen as an extension of Lemma
\ref{odd}.

\begin{lemma}\label{b-lemma} There are constants $a >0, 1>b >0$ such that the following holds.
 Let $W$ be a fixed subspace of
dimension $d \le n-1$ and $X$ a random (row)  vector.  Then

 $$\P( \dist(X,W) \leq
\frac{a}{\sqrt{n}} ) \leq b^{d-n}. $$
\end{lemma}

{\bf \noindent Proof of Lemma \ref{b-lemma}.}   We construct unit
vectors $Z_1, \ldots, Z_{n-d}$ (not necessarily orthogonal) in the
orthogonal complement $W^\perp$ of $W$ as follows.  We let $Z_1$
be an arbitrary unit vector in $W^\perp$; since $Z_1$ has unit
length, at least one of its coordinates has magnitude at least
$1/\sqrt{n}$. Without loss of generality we may assume that it is
the first coordinate $\langle Z_1, e_1 \rangle$ which has
magnitude at least $1/\sqrt{n}$.  Now we let $Z_2$ be an arbitrary
unit vector in $W^\perp \cap e_1^\perp$ (which has dimension at
least $n-d-1$); then $Z_2$ is orthogonal to $e_1$ and has a
coordinate of magnitude at least $1/\sqrt{n}$.  Without loss of
generality we may take $|\langle Z_2, e_2 \rangle| \geq
1/\sqrt{n}$. Continuing in this fashion, we can (without loss of
generality) find $Z_1, \ldots, Z_{n-d} \in W^\perp$ such that each
$Z_j$ is orthogonal to $e_1,\ldots,e_{j-1}$ and is such that
$|\langle Z_j, e_j \rangle| \geq 1/\sqrt{n}$.

Now suppose that $X = (\ep_1,\ldots,\ep_n)$ is such that
$\dist(X,W) \leq \frac{a}{\sqrt{n}}$, where $a$ is a sufficiently
small positive constant. Fix  the last $d$ coordinates
$\ep_{n-d+1},\ldots,\ep_n$ and let $T$ denote the set of all
vectors $X$ with these fixed coordinates satisfying

$$ \dist(X,W) \leq \frac{a}{\sqrt{n}}. $$

 Fix a
vector $X_0 = (g_1, \dots, g_n) $ in $T$. It is easy to show that
for any vector $X =(g_1', \dots, g_n') \in T$, $|g_i'-g_i| \le
2a$, for all $1\le i \le n-d$. On the other hand, if $a$ is
sufficiently small, then by the third property of the $\xi_{ij}$,
there is a positive constant $b <1$ such that the set of $g_i'$
where $|g_i'-g_i| \le 2a$ has measure at most $b$ for all $1 \le i
\le n-d$. This proves the claim. \hs

The rest of the proof is basically the same, with some minor and
natural modification in the calculation.
 The error term obtained  from Lemma \ref{b-lemma} (in the
 determinant) is only

$$n^{- O(n^{1/2} \ln n) } = \exp(-o(n^{1/2+ \ep})) $$

\noindent for any fixed $\ep >0$.

\vskip2mm

In certain situations, we do not have the assumption that
$|\xi_{ij}|$ are bounded from above by a constant. We are going to
consider the following model. Let $\xi_{ij}, 1\le i,j \le n$ be
i.i.d. random variables with mean zero and variance one. Assume
furthermore that their fourth moment is finite. Consider the
random matrix $M_n$ with $\xi_{ij}$ as its entries.

By using Lemmas \ref{odd} and \ref{distance0} and replacing Lemma
\ref{Tal} by a result of Bai and Yin \cite{BY}, which asserts that
the volume of the $(1-\gamma)n$-dimensional parallelepiped spanned
by the first $(1-\gamma)n $ row vectors is at least
$n^{(1/2-\gamma/2 -o(1))n}$ with probability $1-o(1)$ for any
fixed $\gamma >0$, we can prove

\begin{theorem} \label{bai} We have, with probability $1-o(1)$,
that

$$|\det M_n | \ge n^{(1/2-o(1)) n}. $$

\end{theorem}

\vskip2mm

{\it \noindent Acknowledgement.} We would like to thank K. Ball,
J. Bourgain,  N. Linial, K. Maples, A. Naor, G.  Schechtman, G. Ziegler  and
the referees for useful comments.

\end{document}